\long\def\symbolfootnote[#1]#2{\begingroup\def\thefootnote{\fnsymbol{footnote}}\footnote[#1]{#2}\endgroup}

\documentclass{amsart}
\usepackage{amssymb}
\usepackage{amsmath}
\usepackage{amsfonts}

\newtheorem{theorem}{Theorem}[section]
\newtheorem{corollary}[theorem]{Corollary}
\newtheorem{lemma}[theorem]{Lemma}
\theoremstyle{remark}

\theoremstyle{definition}

\theoremstyle{proposition}

\numberwithin{equation}{section}

\begin{document}
\author{Jiaping Wang and Linfeng Zhou}
\title[Eigenform gradient estimate]{Gradient estimate for eigenforms of Hodge Laplacian}
\date{}
\maketitle

\begin{abstract}
In this paper, we derive a gradient estimate for the linear combinations of eigenforms of the Hodge Laplacian on a closed manifold. The estimate is given
in terms of the dimension, volume, diameter and curvature bound of the manifold. As an application, we obtain directly a sharp estimate for the heat kernel of the Hodge Laplacian.
\end{abstract}

\section{Introduction}

\symbolfootnote[0]{The first author is partially supported by an NSF grant.}

Let $(M^n, g)$ be a compact oriented Riemannian manifold without boundary.
The Hodge Laplacian $\Delta: A^p(M)\rightarrow A^p(M),$ acting on
the space of smooth $p$-forms $A^p(M)$ on $M,$ is defined by
\begin{equation*}
\Delta=-d\delta-\delta d.
\end{equation*}
Here, as usual, $d$ is the exterior differential operator and $\delta$ the adjoint of $d$ with respect to the $L^2$ inner product on $A^p(M).$
We denote the eigenvalues of $\Delta$ by
$\{0\leq \lambda_{1}\leq\dots\lambda_{k}\leq\dots\}$ with a corresponding orthonormal basis of eigenforms $\{\phi_i\}_{i=1}^{\infty}.$
We have the following estimate concerning the eigenforms.

\begin{theorem}
Let $(M^n,g)$ be a closed manifold with curvature bound $|Rm|\le K.$
Then for any $b_i\in \mathbb{R}$ with $\sum_{i=1}^k b_i^2\le 1,$ the form
$\omega=\sum_{i=1}^k {b_i}\phi_i$ satisfies the estimate
\begin{equation*}
|\nabla\omega|^2+(\lambda_k+1)\,|\omega|^2\leq c\,(\lambda_k+1)^{\frac{n}{2}+1},
\end{equation*}
where $c=c(n,V,d,K)$ is an explicit constant depending on the dimension $n,$ volume $V,$ diameter $d$ and the curvature bound $K.$
\end{theorem}

We would like to emphasize that the estimate is valid for all finite linear combinations of the eigenforms, and it does not involve any covariant derivatives
of the curvature tensor. Also, the exponent $\frac{n}{2}+1$ in $\lambda_k$ is sharp.
This sharp exponent in turn leads to another one in $k$ for the lower bound of the eigenvalue $\lambda_k\geq c\,k^{-\frac{2}{n}}$ for all $k>b_p,$ the Betti number of the $p-$th cohomology of $M.$

Our estimates can then be applied to analyze the heat kernel of $\Delta.$ Combining
with a result of Rumin \cite{R}, one has the following Sobolev inequality
for $p-$forms.

\begin{theorem} For an explicit constant $c=c(n,V,d,K),$
\begin{equation*}
\left(\int_M |\omega-P(\omega)|^{\frac{2n}{n-2}}\right)^{\frac{n-2}{n}}\le c\, \int_M \{|d \omega|^2+|\delta \omega|^2\}
\end{equation*}
for all smooth $p-$form $\omega$ on $M,$ where $P(\omega)$ denotes the projection
of $\omega$ on to the space of harmonic $p-$forms.
\end{theorem}

Another consequence is the following Hessian estimate for the eigenfunctions on $M.$

\begin{corollary}
Let $(M^n,g)$ be a closed manifold with curvature bound $|Rm|\le K.$
Let $\phi_1,\phi_2,\cdots,\phi_k$ be orthonormal eigenfunctions of the scalar
Laplacian with corresponding eigenvalues $0<\lambda_{1}\le \lambda_{2}\le \cdots \le \lambda_{k}.$ Then there exists a constant $c(K, d, V, n)$ such that
\begin{equation*}
|\nabla d\,f|\leq c \,\lambda_k^{\frac{n+4}{4}},
\end{equation*}
where $f=\sum_{i=1}^k b_i\,\phi_i$ and $\sum_{i=1}^k b_i^2=1.$
\end{corollary}

Let us point out that the analysis of the Laplacian on a compact manifold is a
classical subject. Numerous contributions have been made by various authors. 
While the results here are mostly known, we do hope our seemingly more direct 
treatment is of certain expository value.

As well-known, the gradient estimate method was successfully employed by Yau \cite{Y1} to study harmonic functions on complete manifolds. The method was further developed by Li \cite{L1}, and Li and Yau \cite{L-Y} to study eigenfunctions and eigenvalues. In particular, they have obtained a lower bound for the first non-zero eigenvalue of the scalar Laplacian in terms of the lower bound of Ricci curvature and the diameter of the manifold.

Our current work is very much motivated by and follows the ideas in a famous paper of Li \cite{L2}, where he has obtained a lower bounds for higher eigenvalues of the Hodge Laplacian. This is achieved through an estimate of the linear combinations of
the eigenforms. The estimate involves the curvature operator lower bound
and the Sobolev constant of the manifold, but not the curvature upper bound.
However, the estimate there seems insufficient to provide a sharp exponent for the eigenvalue lower bounds alluded above. We would also like to point out that both E. Aubry's PhD thesis and the paper \cite{B-B-C} by W. Ballmann, J. Br\"uning and G. Carron have already developed a gradient estimate for individual eigenforms. 

The case of scalar Laplacian deserves special attention as it is of more common
concern. So we will treat it separately in section 2. The result is a bit stronger in the sense it only involves the Ricci curvature lower bound in all the estimates.
The approach is also more straightforward as it only relies on a direct application of the maximum principle.

The case of general Hodge Laplacian is handled in section 3. The proof now involves
an iteration scheme as in \cite{L2}.

Finally, we mention that the results here can be extended to the case of compact manifolds with boundary. For the ease of exposition, we omit the details here.

\bigskip

\noindent{\bf Acknowledgements.} We would like to thank Gilles Carron for his insightful comments which lead to various improvement to the paper. Part of the paper was written while the second author was visiting the School of Mathematics at the University of Minnesota. He deeply appreciates its hospitality. He would also like to thank Gang Liu for his helpful comments.

\section{Analysis of scalar Laplacian}

In this section, we will derive a variant version of the well-known gradient estimates of Li-Yau\cite{L-Y} concerning the eigenfunctions. As an application,
we give direct proofs to some well-known results including a lower bound of 
the high eigenvalue, the existence of heat kernel and its long time decay estimate.

Let $(M^n,g)$ be a closed Riemannian manifold with diameter $d,$ volume $V,$
and Ricci curvature lower bound  $-(n-1)K,$ where $K\ge 0$ is a constant.
Denote the eigenvalues of the Laplacian $\Delta$ by $0=\lambda_0<\lambda_1\leq\cdots\leq\lambda_k\leq \cdots $ with the corresponding eigenfunction $\phi_i,$ $i=0,1,2,\cdots,$ satisfying
\begin{equation*}
\Delta \phi_i=-\lambda_i\phi_i,\quad \int_{M}\phi_i \,\phi_j=\delta_{ij}.
\end{equation*}

For a given constant $c,$ consider the function
\begin{equation*}
Q(x)=|\nabla \phi|^2+c\,\phi^2,
\end{equation*}
where $\phi=\sum_{i=1}^k b_i\phi_i$ with $b_i\in \mathbb{R}$ and $\sum_{i=1}^k b_i^2=1.$
Obviously, the maximum value of $Q(x)$ over $M$ is a function of $b_1,\cdots,b_k.$
This function in turn achieves its maximum at some point $a_1,\cdots,a_k.$ Let $u=\sum_{i=1}^k a_i\phi_i.$

\begin{lemma}\label{le1}
\begin{equation*}
|\nabla u|^2+A\,u^2\leq A\, \max_{M}u^2,
\end{equation*}
where $A=\lambda_k+(n-1)K$.
\end{lemma}

\begin{proof} Define
\begin{equation*}
F(b_1,\dots,b_k,x,\lambda)=Q(x)-\lambda \,(\sum_{i=1}^k b_i^2-1).
\end{equation*}
Then, subject to the constraint $\sum_{i=1}^k b_i^2=1,$ $F$ achieves its maximum value at some point $(a_1,\cdots,a_k,x_0,\alpha).$
We now show
\begin{equation*}
|\nabla u|^2(x_0)+c\,u^2(x_0)\leq c\,\max_{M} u^2
\end{equation*}
for $c> \lambda_k+(n-1)K.$

At the point  $(a_1,\cdots,a_k,x_0,\alpha)$, $F$ satisfies

\begin{equation}\label{eq2}\left\{\begin{array}{l}
\nabla F(a_1,\cdots,a_k,x_0,\alpha)=0 \\
\Delta F(a_1,\cdots,a_k,x_0,\alpha)\leq 0\\
\frac{\partial F}{\partial b_i}=0\\
\sum_{i=1}^k a_i^2=1.
\end{array}\right.\end{equation}
From the third equation of (\ref{eq2}), we have

\begin{equation*}
\sum_{j=1}^{k} (2a_j\langle\nabla \phi_i,\nabla \phi_j\rangle+2ca_j\langle\phi_i,\phi_j\rangle)-2\alpha a_i=0.
\end{equation*}
After multiplying by $a_i$ and summing over $i,$ one sees that

\begin{equation*}
\alpha=Q(u,x_0)=|\nabla u|^2(x_0)+c\,u^2(x_0).
\end{equation*}

Suppose now that
\begin{equation*}
|\nabla u|^2(x_0)+cu^2(x_0)>c\,\max_{M}\,u^2.
\end{equation*}
Then
\begin{equation*}
\nabla u(x_0)\neq 0
\end{equation*}
and one can choose an orthonormal frame $\{e_1,\dots,e_n\}$ at $x_0$ so that
\begin{equation*}
\nabla u(x_0)=u_1(x_0)e_1.
\end{equation*}
Now the first equation of (\ref{eq2}), $\nabla F(a_1,\dots,a_k,x_0,\alpha)=0,$
becomes
\begin{equation*}
2u_{1}u_{1i}+2cuu_{i}=0
\end{equation*}
for $i=1,\dots,n$.  This in particular implies

\begin{equation}\label{eq3}
|\nabla\nabla u|^2\ge u_{11}^2=c^2\,u^2.
\end{equation}

On the other hand, at the maximum point $(a_1,\cdots,a_k,x_0,\alpha),$
\begin{equation*}
\Delta F(a_1,\dots,a_k,x_0,\alpha)\leq  0
\end{equation*}
or
\begin{equation}\label{eq4}
\Delta |\nabla u|^2+c\Delta u^2\leq 0.
\end{equation}
By the Bochner formula, it becomes
\begin{equation*}
|\nabla\nabla u|^2+\langle \nabla\Delta  u, \nabla u\rangle+\langle Ric(\nabla u,\nabla u)\rangle+cu\Delta u+c|\nabla u|^2 \le 0.
\end{equation*}
In view of (\ref{eq3}) and the lower bound of Ricci curvature, the above inequality
reduces to
\begin{equation*}
c^2u^2+\langle \nabla\Delta u,\nabla u\rangle-(n-1)K|\nabla u|^2+cu\Delta u+c|\nabla u|^2\le 0.
\end{equation*}
Note that
\begin{equation*}
\Delta u=-\sum_{i=1}^k \lambda_i \,a_i\,\phi_i.
\end{equation*}
Therefore,
\begin{eqnarray*}
0&\geq& c^2\,u^2+(c-(n-1)K)|\nabla u|^2-\sum_{i,j=1}^k\lambda_i\, a_i\, a_j\,\langle\nabla\phi_i,\nabla\phi_j\rangle-\sum_{i,j=1}^k c\,\lambda_i\,a_i\,a_j\,\phi_i,\phi_j\\
&\geq& c^2\,u^2+(c-(n-1)K)|\nabla u|^2
-\sum_{i=1}^k\,\lambda_i\,a_i\,\sum_{j=1}^k(a_j\langle\nabla\phi_i,\nabla\phi_j\rangle
+c\,a_j\,\phi_i\,\phi_j)\\
&\geq& c^2\,u^2+(c-(n-1)K)|\nabla u|^2-\sum_{i=1}^k\,\alpha\,\lambda_i\,a_i^2\\
&\geq& c^2\,u^2+(c-(n-1)K)|\nabla u|^2-\alpha\,\lambda_k\\
&\geq& c^2\,u^2+(c-(n-1)K)|\nabla u|^2-\lambda_k\,(|\nabla u|^2+cu^2)\\
&\geq& c\,(c-\lambda_k)\,u^2+(c-(n-1)K-\lambda_k)|\nabla u|^2.
\end{eqnarray*}
This is an obvious contradiction if $c>(n-1)K+\lambda_k.$
In other words,
\begin{equation*}
|\nabla u|^2(x_0)+c\,u^2(x_0)\leq c\,\max_{M}\,u^2
\end{equation*}
for all $c>(n-1)K+\lambda_k.$
The lemma follows by letting $c$ approach $\lambda_k+(n-1)K.$
\end{proof}

As a consequence, we obtain a quick proof to the following well-known facts.

\begin{theorem}\label{th1}
There exists a constant $c(K, d, V, n)$ such that

\noindent (1)
\begin{equation*}
|\nabla \phi|^2\leq c \lambda_k^{\frac{n+2}{2}}, \quad \phi^2\leq c\lambda_k^{\frac{n}{2}}.
\end{equation*}
In particular,
\begin{equation*}
|\nabla \phi_k|\leq c \lambda_k^{\frac{n+2}{4}}, \quad |\phi_k|\leq c\lambda_k^{\frac{n}{4}}.
\end{equation*}

\noindent (2) For all $k\ge 1,$
\begin{equation*}
\lambda_k\geq c^{-1}\,k^{\frac{2}{n}}.
\end{equation*}

\noindent (3) The function $H(x,y,t)$ given by

\begin{equation*}
H(x,y,t)=\frac{1}{V}+\sum_{k=1}^\infty e^{-\lambda_k t}\,\phi_k(x) \,\phi_k(y)
\end{equation*}
is a heat kernel of $M.$ Moreover,
\begin{equation*}
|H(x,y,t)-\frac{1}{V}|\le c\, t^{-\frac{n}{2}}
\end{equation*}
for all $t>0.$

\noindent (4) The following Sobolev inequality holds.

\begin{equation*}
\left(\int_M |f|^{\frac{2n}{n-2}}\right)^{\frac{n-2}{n}}\le c\, \int_M |\nabla f|^2
\end{equation*}
for all smooth function $f$ on $M$ with $\int_M f=0.$
\end{theorem}

\begin{proof} (1) Let $u$ be the function considered in the preceding lemma. Then we need only to prove the estimate for $u.$ Choose point $p$ such that

\begin{equation*}
u^2(p)=\max_{M} \,u^2.
\end{equation*}
For $r>0$ and $x\in B_p(\frac{r}{\sqrt{\lambda_k+(n-1)K}}),$ using lemma \ref{le1}, we conclude

\begin{eqnarray*}
u^2(p)-u^2(x)&\leq& \max_{y\in M} 2|u|(y)\,|\nabla u|(y)\,d(x,p)\\
&\leq& 2\,u^2(p)\,\sqrt{\lambda_k+(n-1)K}\,\frac{r}{\sqrt{\lambda_k+(n-1)K}}\\
&\leq& 2r\,u^2(p).
\end{eqnarray*}
Therefore,

\begin{equation*}
u^2(x)\geq(1-2r)\,u^2(p)
\end{equation*}
on $B_p(\frac{r}{\sqrt{\lambda_k+(n-1)K}}).$
Integrating with respect to $x$ over the ball yields

\begin{equation*}
1=||u||^2_{L^2(M)}\geq (1-2r)\,u^2(p)\,\frac{V_p(\frac{r}{\sqrt{\lambda_k+(n-1)K}})}{V_p(d)}\,V_p(d).
\end{equation*}
Choose $r$ such that

\begin{equation*}
4r<1 \quad\text{and}\quad \frac{r}{\lambda_1+\sqrt{(n-1)K}}<d.
\end{equation*}
Then by the Bishop volume comparison theorem we have

\begin{equation*}
1\geq(1-2r)\, u^2(p)\,\frac{c(K,d,V,n)}{(\lambda_k+(n-1)K)^{\frac{n}{2}}}.
\end{equation*}
In other words,

\begin{equation*}
u^2(x)\leq u^2(p) \leq c(K,d,V,n)\,\lambda_k^{\frac{n}{2}},
\end{equation*}
where we have used the fact that $\lambda_k\ge \lambda_1\ge c$ by \cite{L-Y}.
Using lemma \ref{le1} again, we also conclude

\begin{equation*}
|\nabla u|^2(x) \leq (\lambda_k+(n-1)K)\, u^2(p) \leq c(K,d,V,n)\,\lambda_k^{\frac{n+2}{2}}.
\end{equation*}

(2) For each $x\in M,$ there exists an orthogonal matrix
$(a_{ij})_{k\times k}$ such that

\begin{equation*}
\nabla \psi_j(x)=0
\end{equation*}
for $j=n+1,\cdots,k,$ where $\psi_j=\sum_{i=1}^k a_{ij}\,\phi_i.$

From (1), it follows that
\begin{eqnarray*}
\sum_{i=1}^k|\nabla \phi_i|^2(x)&=&\sum_{j=1}^n |\nabla \psi_j|^2(x)\\
&\leq& n\,\max_j\,|\nabla \psi_j|^2 \\
&\leq& c_1\,\lambda_k^{\frac{n+2}{2}}.
\end{eqnarray*}
Integrating the inequality with respect to $x,$ we conclude

\begin{equation*}
\lambda_1+\lambda_2+\dots+\lambda_k\leq c_2\,\lambda_k^{\frac{n+2}{2}}.
\end{equation*}
By an elementary induction argument, the inequality implies
\begin{equation*}
\lambda_k\geq c_3\,k^{2/n}
\end{equation*}
for all $k\geq 1,$ where $c_3=\min\{\lambda_1, (\frac{1}{c_2}\frac{n}{n+2})^\frac{n}{2} \}.$

(3) In view of (1) and (2), it is straightforward to check the infinite
series

\begin{equation*}
\frac{1}{V}+\sum_{k=1}^\infty e^{-\lambda_k t}\,\phi_k(x)\, \phi_k(y)
\end{equation*}
converges uniformly in the $C^1$ sense for $x, y\in M$ and $t\ge c$
for any $c>0.$ It is then easy to verify the limit is a heat kernel of $M.$

Since
\begin{equation*}
c^{-1}\,k^{\frac{2}{n}}\le \lambda_k\leq c\,k^{\frac{2}{n}},
\end{equation*}
one sees by (1) that
\begin{eqnarray*}
|H(x,y,t)-\frac{1}{V}|&\le& \sum_{k=1}^\infty e^{-\lambda_k t}|\phi_k|(x)\, |\phi_k|(y)\\
&\le&  \sum_{k=1}^\infty c\, \lambda_k^{\frac{n}{2}}\,e^{-\lambda_k t}\\
&\le& c\,t^{-\frac{n}{2}}\,\int_0^\infty s^{\frac{n}{2}}\,e^{-s}\,ds\\
&\le& c\, t^{-\frac{n}{2}}.
\end{eqnarray*}

(4) follows from (3) by a result of Varopoulos \cite{V}.
\end{proof}

We remark that both (2) and (3) were first proved by Cheng and Li \cite{C-L}
using the Sobolev inequality. Historically, the Sobolev inequality on manifolds
was derived from the isoperimetric inequalities, which were established
by Yau \cite{Y2} and Croke \cite{C}.

\section{Gradient estimate for eigenforms}
Using the well-known Bochner-Weitzenbock formula, one can directly apply
the proof in the previous section to the Hodge Laplacian acting on the
smooth $p-$forms on $M.$ However, the resulting estimates depend also
on the bounds of the covariant derivative of the curvature tensor of $M.$
It turns out this dependency is superfluous by adopting a different
argument as demonstrated by E. Aubry in his PhD thesis and also by W. Ballmann, J. Br\"uning
and G. Carron in \cite{B-B-C}. In the following,
we present a slightly refined version of their argument to suit our purpose.

We will use the moving frame notations. So for a $p-$form
$\omega$ on $M,$ under an orthonormal coframe $\{\omega_1,\cdots,\omega_n\},$  $\omega=a_{i_1\cdots i_p}\omega_{i_p}\wedge\cdots \wedge\omega_{i_1}.$

The Bochner-Weitzenbock formula says
\begin{equation*}
\Delta \omega=\Delta_B \omega-E(\omega),
\end{equation*}
where
\begin{equation*}
\Delta_B \omega=a_{i_1\cdots i_p,jj}\omega_{i_p}\wedge\cdots \wedge\omega_{i_1}
\end{equation*}
is the Bochner Laplacian and
\begin{equation*}
E(\omega)=
R_{k_{\beta}i_{\beta}j_{\alpha}i_{\alpha}}\,a_{i_1\dots k_{\beta}\dots i_p}\,\omega_{i_p}\wedge\cdots \wedge\omega_{j_{\alpha}}\wedge\cdots\wedge\omega_{i_1}
\end{equation*}
with $R_{ijkl}$ being the curvature tensor of $M.$
Now,
\begin{equation*}
\Delta_B (\nabla \omega)=a_{i_1\dots i_p,ijj}\omega_{i_p}\wedge\cdots \wedge\omega_{i_1}\otimes \omega_i
\end{equation*}
and
\begin{equation*}
\nabla \Delta \omega=a_{i_1\dots i_p,jji}\omega_{i_p}\wedge\cdots \wedge\omega_{i_1}\otimes \omega_i-\nabla (E(\omega)).
\end{equation*}
Hence
\begin{equation*}
\Delta_B (\nabla \omega)-\nabla\Delta \omega=a_{I,ijj}\omega_{I}\otimes \omega_i-a_{I,jji} \omega_{I}\otimes \omega_i+\nabla (E(\omega)).
\end{equation*}
By the Ricci identity, we have
\begin{equation*}
a_{I,ijj}-a_{I,jij}=(R_{j_\alpha i_\alpha ij}a_{i_1\dots j_\alpha\dots i_p})_{,j}
\end{equation*}
and
\begin{equation*}
a_{I,jij}-a_{I,jji}=R_{j_\alpha i_\alpha ij}a_{i_1\dots j_\alpha \dots i_p,j}+R_{ljij}a_{i_1\dots i_p,l}.
\end{equation*}
Thus we have the commutation formula
\begin{eqnarray}
\Delta_B(\nabla \omega)-\nabla\Delta \omega&=&R_{li}a_{i_1\dots i_p,l}\,\omega_{i_p}\wedge\cdots \wedge\omega_{i_1}\otimes \omega_i \label{eq5} \\
&+& R_{j_\alpha i_\alpha ij}a_{i_1\dots j_\alpha \dots i_p,j}\,\omega_{i_p}\wedge\cdots \wedge\omega_{i_1}\otimes \omega_i \notag \\
&+& \left(R_{j_\alpha i_\alpha ij}\,a_{i_1\dots j_\alpha \dots i_p}\right)_{,j}\,\omega_{i_p}\wedge\dots \wedge\omega_{i_1}\otimes \omega_i \notag\\
&+& \nabla (E(\omega)). \notag
\end{eqnarray}
Finally, we conclude
\begin{eqnarray}
\langle\Delta_B(\nabla\omega)-\nabla\Delta\omega,\nabla\omega\rangle
&=& R_{li}a_{i_1\dots i_p,l}a_{i_1\dots i_p,i} \label{eq6}\\
&+& R_{j_\alpha i_\alpha ij}a_{i_1\dots j_\alpha \dots i_p,j}a_{i_1\dots i_\alpha\dots i_p,i} \notag\\
&+& \left(R_{j_\alpha i_\alpha ij}\,a_{i_1\dots j_\alpha \dots i_p}\right)_{,j}a_{i_1\dots i_\alpha\dots i_p,i} \notag \\
&+& \langle \nabla (E(\omega)),\nabla\omega\rangle \notag
\end{eqnarray}
Note that these formulas and the following lemma have more or less been derived by Le Couturier and G. Robert in \cite{L-R}.

We now consider the function $f=|\nabla\omega|^2+A\,|\omega|^2,$
where $A\geq 1$ is a fixed constant.

\begin{lemma}\label{le6}
Let $(M^n,g)$ be a closed Riemannian manifold with curvature operator
$|Rm|\leq K.$ Then for $k\geq 1,$

\begin{equation*}
\int_M f^{k-1}\Delta f \geq 2\,\int_M\left(\langle \nabla\Delta\omega, \nabla\omega\rangle+A\,\langle\Delta\omega,\omega\rangle\right)\,f^{k-1}
- c\, k^2\,\int_M f^{k},
\end{equation*}
where $c=2nK(K+2)+18\,K^2.$
\end{lemma}

\begin{proof}
Direct calculation gives
\begin{eqnarray}
\Delta f&=&\Delta \left(|\nabla\omega|^2+A\,|\omega|^2\right) \label{eq7}\\
&=& 2\langle\Delta_B(\nabla\omega), \nabla\omega\rangle+2|\nabla\nabla\omega|^2 \notag\\
&+&2A\,|\nabla\omega|^2+2A\,\langle \Delta_B \omega,\omega\rangle \notag \\
&=& 2\langle \nabla \Delta \omega, \nabla\omega\rangle+2A\,\langle \Delta \omega,\omega\rangle \notag\\
&+& 2|\nabla\nabla\omega|^2+ 2A\,|\nabla\omega|^2+2A\,\langle E(\omega),\omega\rangle \notag \\
&+& 2\,\langle\Delta_B(\nabla\omega)-\nabla \Delta \omega, \nabla\omega\rangle. \notag
\end{eqnarray}
Therefore,
\begin{eqnarray}
\int_{M}f^{k-1}\Delta f&=& 2\,\int_{M}\left(\langle\nabla\Delta\omega,\nabla\omega\rangle
+A\,\langle\Delta\omega,\omega\rangle\right)f^{k-1} \label{eq8}\\
&+& 2\,\int_M \left(|\nabla\nabla\omega|^2+ A\,|\nabla\omega|^2\right)\,f^{k-1} \notag\\
&+& 2A\,\int_M \langle E(\omega),\omega\rangle \, f^{k-1} \notag\\
&+& 2\,\int_M \langle\Delta_B(\nabla\omega)-\nabla \Delta \omega, \nabla\omega\rangle \,f^{k-1}.\notag
\end{eqnarray}
Since $|Rm|\leq K,$

\begin{equation}
2A\,\int_M \langle E(\omega),\omega\rangle \, f^{k-1}\geq
-2K\,\int_M f^k.\label{eq9}
\end{equation}
Using (\ref{eq6}), we have
\begin{eqnarray}
&&2\,\int_M\langle\Delta_B(\nabla\omega)-\nabla \Delta \omega, \nabla\omega\rangle \,f^{k-1} \label{eq10}\\
&=& 2\,\int_M R_{li}a_{i_1\dots i_p,l}a_{i_1\dots i_p,i}\,f^{k-1} \notag\\
&+& 2\,\int_M R_{j_\alpha i_\alpha ij}a_{i_1\dots j_\alpha \dots i_p,j}a_{i_1\dots i_\alpha\dots i_p,i}\,f^{k-1} \notag\\
&+& 2\,\int_M \left(R_{j_\alpha i_\alpha ij}\,a_{i_1\dots j_\alpha \dots i_p}\right)_{,j}a_{i_1\dots i_\alpha\dots i_p,i}\,f^{k-1} \notag\\
&+& 2\,\int_M \langle \nabla (E(\omega)),\nabla\omega\rangle f^{k-1}. \notag
\end{eqnarray}
The first and second term of (\ref{eq10}) obviously satisfy
\begin{equation}
2\,\int_M R_{li}a_{i_1\dots i_p,l}a_{i_1\dots i_p,i}\,f^{k-1}\geq
-2(n-1)K\,\int_M f^k.\label{eq11}
\end{equation}
and
\begin{eqnarray}
2\,\int_M R_{j_\alpha i_\alpha ij}a_{i_1\dots j_\alpha \dots i_p,j}a_{i_1\dots i_\alpha\dots i_p,i}\,f^{k-1} &\geq& -2K\,\int_M |\nabla \omega|^2\,f^{k-1}
\label{eq12} \\
&\geq& -2K\,\int_M f^k. \notag
\end{eqnarray}
For the third term of (\ref{eq10}), after integration by parts, we have
\begin{eqnarray}
&&2\,\int_M \left(R_{j_\alpha i_\alpha ij}\,a_{i_1\dots j_\alpha \dots i_p}\right)_{,j}a_{i_1\dots i_\alpha\dots i_p,i}\,f^{k-1} \label{eq13}\\
&=&-2\,\int_M R_{j_\alpha i_\alpha ij}\,a_{i_1\dots j_\alpha \dots i_p}
a_{i_1\dots i_\alpha\dots i_p,ij}\,f^{k-1} \notag \\
&-&2(k-1)\,\int_M R_{j_\alpha i_\alpha ij}\,a_{i_1\dots j_\alpha \dots i_p}
a_{i_1\dots i_\alpha\dots i_p,i}\,f^{k-2}\,f_j \notag\\
&\geq& -2K\,\int_M |\omega|\,|\nabla \nabla \omega| \,f^{k-1} \notag\\
&-&2(k-1)K\,\int_M |\omega|\,|\nabla \omega| \,f^{k-2}\,|\nabla f| \notag\\
&\geq& -2K^2\,\int_M f^k-\frac{1}{2}\,\int_M |\nabla \nabla \omega|^2 \,f^{k-1} \notag\\
&-& 8\,k^2\,K^2\,\int_M f^{k}-\frac{1}{2}\,\int_M (|\nabla \nabla \omega|^2+A\,|\nabla \omega|^2) \,f^{k-1}, \notag
\end{eqnarray}
where we have used the fact that
\begin{equation*}
|\omega|\,|\nabla \omega|\leq f
\end{equation*}
and
\begin{eqnarray}
|\nabla f|&\leq& 2|\nabla \omega|\, |\nabla \nabla \omega|+2A\,|\omega|\, |\nabla \omega| \label{eq14}\\
&\leq& 4k\,K\,(|\nabla \omega|^2+A\, |\omega|^2)+\frac{1}{4k\,K}(|\nabla \nabla \omega|^2+A\,|\nabla \omega|^2). \notag
\end{eqnarray}
Applying integration of parts to the last term of (\ref{eq10}), we get
\begin{eqnarray}
&&2\,\int_M \langle \nabla (E(\omega)),\nabla\omega\rangle f^{k-1} \label{eq15}\\
&\geq& -2\,\int_M \langle E(\omega),\Delta_B \omega\rangle f^{k-1} \notag \\
&-& 2(k-1)\,\int_M |E(\omega)|\,|\nabla\omega|\, f^{k-2}\,|\nabla f| \notag\\
&\geq&-2\sqrt{n}\,K\,\int_M |\omega|\,|\nabla \nabla \omega|\, f^{k-1} \notag\\
&-& 2(k-1)K\,\int_M |\omega|\,|\nabla\omega|\, f^{k-2}\,|\nabla f| \notag\\
&\geq&-2n\,K^2\,\int_M f^{k}- \frac{1}{2}\,\int_M |\nabla \nabla \omega|^2\,f^{k-1}\notag\\
&-& 8\,k^2\,K^2\,\int_M f^{k}-\frac{1}{2}\,\int_M (|\nabla \nabla \omega|^2+A\,|\nabla \omega|^2) \,f^{k-1}, \notag
\end{eqnarray}
where we have used (\ref{eq14}) in the last step.

Putting (\ref{eq11}), (\ref{eq12}), (\ref{eq13}) and (\ref{eq15}) into (\ref{eq10}),
we conclude

\begin{eqnarray}
&&2\,\int_M \langle\Delta_B(\nabla\omega)-\nabla \Delta \omega, \nabla\omega\rangle \,f^{k-1}  \label{eq16}\\
&\geq& -\left(2nK(K+1)+18\,k^2\,K^2\right)\,\int_M f^k \notag\\
&-& 2\,\int_M (|\nabla \nabla \omega|^2+A\,|\nabla \omega|^2) \,f^{k-1}. \notag
\end{eqnarray}
Plugging (\ref{eq9}) and (\ref{eq16}) into (\ref{eq8}), we arrived at
\begin{eqnarray*}
\int_{M}f^{k-1}\Delta f &\geq& 2\,\int_{M}\left(\langle\nabla\Delta\omega,\nabla\omega\rangle
+A\,\langle\Delta\omega,\omega\rangle\right)\,f^{k-1}\\
&-&\left(2nK(K+2)+18\,k^2\,K^2\right)\,\int_M f^k.
\end{eqnarray*}
The lemma is proved.
\end{proof}

We now prove a gradient estimate concerning the linear combinations of eigenforms.

\begin{theorem} \label{th4}
Let $(M^n,g)$ be a closed manifold with curvature bound $|Rm|\le K.$
Let $\phi_1,\phi_2,\cdots,\phi_l$ be orthonormal eigenforms of the Hodge Laplacian $\Delta$ acting on the $p-$forms with corresponding eigenvalues $0\le \lambda_{1}\le\lambda_{2}\le \cdots \le \lambda_{l}.$
Then for any $b_i\in \mathbb{R}$ with $\sum_{i=1}^l b_i^2\le 1,$ the form
$\omega=\sum_{i=1}^l{b_i}\phi_i$ satisfies the estimate
\begin{equation*}
|\nabla\omega|^2+A\,|\omega|^2\leq c\,(\lambda_l+K+1)^{\frac{n}{2}+1},
\end{equation*}
where  $A=\lambda_l+K+1,$ and $c=c(n,V,d,K)$ is a constant.
\end{theorem}

\begin{proof} For each $k\geq 1,$ let

\begin{equation*}
I_k=\max \int_M f^{2k},
\end{equation*}
where $f=|\nabla \omega|^2+A\,|\omega|^2$ and
the maximum is taken over all $\omega=\sum_{i=1}^l b_i\,\phi_i$ with $b_i\in \mathbb{R}$ and $\sum_{i=1}^l b_i^2\le 1.$

Note that for $\omega=\sum_{i=1}^l{b_i}\phi_i$ with $\sum_{i=1}^l b_i^2\le 1,$
\begin{equation*}
\Delta \omega=-\sum_{i=1}^l \lambda_i\,b_i \,\phi_i=-\lambda_l \,\sum_{i=1}^l{a_i}\phi_i,
\end{equation*}
where $a_i=\lambda_i\,\lambda_l^{-1}\,b_i,$ $i=1,\cdots, l.$
Obviously,

$$
\sum_{i=1}^l a_i^2\le 1.
$$
So if we denote $\eta=\sum_{i=1}^l a_i\,\phi_i,$ then

\begin{eqnarray*}
&&\int_M \left(\langle \nabla\Delta\omega, \nabla\omega\rangle+A\,\langle\Delta\omega,\omega\rangle\right)f^{2k-1}\\
&\geq& -\int_M \left(|\nabla\Delta\omega|^2+A\,|\Delta\omega|^2\right)^{\frac{1}{2}}\,
f^{2k-\frac{1}{2}}\\
&\ge& -\lambda_l\,\left(\int_M (|\nabla\eta|^2+A\,|\eta|^2)^{2k}\right)^{\frac{1}{4k}}\,
\left(\int_M f^{2k}\right)^{\frac{4k-1}{4k}}\\
&\ge& -\lambda_l\,I_k.
\end{eqnarray*}
Combining with lemma \ref{le6}, we have the estimate
\begin{equation}
\int_{M}f^{2k-1}\Delta f\geq -(2\,\lambda_l+c_1\,k^2)\,I_k, \label{eq20}
\end{equation}
where $c_1=8nK(K+2)+72\,K^2.$

On the other hand
\begin{equation}
\int_M f^{2k-1}\Delta f=-\frac{2k-1}{k^2}\int_M |\nabla f^{k}|^2.\label{eq21}
\end{equation}
Applying the Sobolev inequality
\begin{equation*}
\left(\int_M |u|^{2\beta}\right)^{\frac{1}{\beta}}\le C_s\,\left(\int_M |\nabla u|^2
+\int_M u^2\right),
\end{equation*}
where $\beta=\frac{n}{n-2},$ to $u=f^k,$ we get
\begin{equation}
\left(\int_M f^{2k\beta}\right)^{\frac{1}{\beta}}\le C_s\,\left(\int_M |\nabla f^k|^2+\int_M f^{2k}\right). \label{eq22}
\end{equation}
Combining (\ref{eq20}), (\ref{eq21}) and (\ref{eq22}), we get
\begin{equation*}
\left(\int_M f^{2k\beta}\right)^{\frac{1}{\beta}}\leq C_s\,k\,(\lambda_l+c_1\,k^2)\,I_k
\end{equation*}
Since this is true for all $\omega,$ we may maximize the left hand side over $\omega$ and conclude
\begin{equation*}
\left(I_{\beta k}\right)^{\frac{1}{\beta k}}\leq
\left(C_s\,k\,(\lambda_l+c_1\,k^2)\right)^{\frac{1}{k}}\,\left(I_k\right)^{\frac{1}{k}}
\end{equation*}
for all $k\geq 1.$

Let $k=\beta^i,$ $i=0,1,2,\cdots$ and iterate the preceding inequality.
Then,
\begin{eqnarray*}
\lim_{i\to \infty}\left(I_{\beta^i}\right)^{\frac{1}{\beta^i}}
&\leq& \prod_{i=0}^{\infty}\left(C_s\,\beta^i\,(\lambda_l+c_1\,\beta^{2i})\right)^{\frac{1}{\beta^i}}
\,I_{1}\\
&\leq& c_2\,(\lambda_l+1)^{\frac{n}{2}}\,I_{1},
\end{eqnarray*}
where $c_2=c_2(n,d,V,K)$ is a constant. In other words,
\begin{eqnarray*}
&&\max_{\omega}\,\max_{x\in M} \left(|\nabla \omega|^2+A\,|\omega|^2\right)^2(x)\\
&\leq& c_2\,(\lambda_l+1)^{\frac{n}{2}}\,
\max_{\omega} \int_M (|\nabla \omega|^2+A\,|\omega|^2)^2\\
&\leq& c_2\,(\lambda_l+1)^{\frac{n}{2}}\,\max_{\omega}\,\max_{x\in M} \left(|\nabla \omega|^2+A\,|\omega|^2\right)(x)\,\max_{\omega} \int_M (|\nabla \omega|^2+A\,|\omega|^2).
\end{eqnarray*}
Hence,
\begin{equation*}
\max_{\omega}\,\max_{x\in M} \left(|\nabla \omega|^2+A\,|\omega|^2\right)(x)
\leq c_2\,(\lambda_l+1)^{\frac{n}{2}}\,\max_{\omega} \int_M (|\nabla \omega|^2+A\,|\omega|^2).
\end{equation*}
However,
\begin{eqnarray*}
&&\int_M \left(|\nabla\omega|^2+A\,|\omega|^2\right)\\
&=&-\int_M\langle\Delta\omega,\omega\rangle-\int_M\langle E(\omega),\omega\rangle
+c\,\int_M |\omega|^2\\
&\leq&\lambda_l+K+A.
\end{eqnarray*}
Therefore,
\begin{equation*}
\max_{\omega}\,\max_{x\in M} \left(|\nabla \omega|^2+A\,|\omega|^2\right)(x)\leq
c_2\,(\lambda_l+1)^{\frac{n}{2}}\,(\lambda_l+K+A).
\end{equation*}
The theorem is proved.
\end{proof}

As in section 2, we can draw the following conclusions from theorem \ref{th4}.

\begin{theorem}
Let $(M^n,g)$ be a closed manifold with curvature bound $|Rm|\le K.$
Let $0\le \lambda_{1}\le \lambda_{2}\le \cdots \le \lambda_{k}\le \cdots $ be all the eigenvalues of the Hodge Laplacian $\Delta$ acting on the $p-$forms, and
$\phi_1,\phi_2,\cdots,\phi_k,\cdots $ the corresponding orthonormal eigenforms.
Then there exists a constant $c(K, d, V, n)$ such that

\noindent (1) $|\nabla \phi_k|\leq c \,(\lambda_k+1)^{\frac{n+2}{4}}$ and
 $|\phi_k|\leq c\,(\lambda_k+1)^{\frac{n}{4}}.$

\noindent (2) For all $k>b_p,$ the Betti number of the $p-$th cohomology of $M,$
\begin{equation*}
\lambda_k\geq c^{-1}\,k^{\frac{2}{n}}.
\end{equation*}

\noindent (3) The tensor $H_p(x,y,t)$ given by

\begin{equation*}
H_p(x,y,t)=\sum_{k=1}^\infty e^{-\lambda_k t}\phi_k(x) \otimes \phi_k(y)
\end{equation*}
is a heat kernel of $\Delta.$ Moreover,
\begin{equation*}
|H_p(x,y,t)-\sum_{k=1}^{b(p)} \phi_k(x)\otimes \phi_k(y)|\le c\, t^{-\frac{n}{2}}
\end{equation*}
for all $t>0.$

\noindent (4) The following Sobolev inequality holds.

\begin{equation*}
\left(\int_M |\omega-P(\omega)|^{\frac{2n}{n-2}}\right)^{\frac{n-2}{n}}\le c\, \int_M \{|d \omega|^2+|\delta \omega|^2\}
\end{equation*}
for all smooth $p-$form $\omega$ on $M,$ where $P(\omega)$ denotes the projection
of $\omega$ on to the space of harmonic $p-$forms.
\end{theorem}

\begin{proof} (1) is obvious by theorem \ref{th4}. Using theorem \ref{th4}, (2) follows as in the proof of (2) in theorem \ref{th1}, where we now use a result
of T. Mantuano \cite{M} that $\lambda_{b_p+1}\ge c(n,V,d,K).$
For (3), the proof is the same as (3) in theorem \ref{th1}. Finally,
(4) follows from (3) by Theorem 1.2 in \cite{R}.
\end{proof}

We also have the following corollary concerning the eigenfunctions.

\begin{corollary}
Let $(M^n,g)$ be a closed manifold with curvature bound $|Rm|\le K.$
Let $\phi_1,\phi_2,\cdots,\phi_k$ be orthonormal eigenfunctions of the scalar
Laplacian with corresponding eigenvalues $0<\lambda_{1}\le \lambda_{2}\le \cdots \le \lambda_{k}.$ Then there exists a constant $c(K, d, V, n)$ such that
\begin{equation*}
|\nabla d\,\phi|\leq c \,\lambda_k^{\frac{n+4}{4}},
\end{equation*}
where $\phi=\sum_{i=1}^k b_i\,\phi_i$ and $\sum_{i=1}^k b_i^2=1.$
\end{corollary}

\begin{proof} Note that $d\phi_i$ is an eigenform for the Hodge Laplacian
acting on the one forms. Now the corollary follows by applying theorem \ref{th4}
to the one form setting with $d\phi_i$ normalized to have unit length in the $L^2$
sense.
\end{proof}

As a final remark, it is possible to make explicit of all the constants in our arguments. In particular, we could spell out their dependency on the geometric quantities $d,$ $V$ and $K.$

{\small SCHOOL OF MATHEMATICS, UNIVERSITY OF MINNESOTA }

{\small MINNEAPOLIS, MN 55455} \newline
{\small E-mail address: jiaping@math.umn.edu}

\bigskip

{\small DEPARTMENT OF MATHEMATICS, EAST CHINA NORMAL UNIVERSITY }

{\small SHANGHAI 200062, CHINA}\newline
{\small E-mail address: lfzhou@math.ecnu.edu.cn}
\end{document}